\documentclass[12pt]{article}
\usepackage{graphicx,psfrag,amsmath,amssymb,latexsym,color,subfigure,hyperref}

\newcommand{\BEAS}{\begin{eqnarray*}}
\newcommand{\EEAS}{\end{eqnarray*}}
\newcommand{\BEA}{\begin{eqnarray}}
\newcommand{\EEA}{\end{eqnarray}}
\newcommand{\BEQ}{\begin{equation}}
\newcommand{\EEQ}{\end{equation}}
\newcommand{\BIT}{\begin{itemize}}
\newcommand{\EIT}{\end{itemize}}
\newcommand{\BNUM}{\begin{enumerate}}
\newcommand{\ENUM}{\end{enumerate}}

\newcommand{\BA}{\begin{array}}
\newcommand{\EA}{\end{array}}

\newcommand{\ones}{\mathbf 1}

\newcommand{\reals}{{\mbox{\bf R}}}

\newcommand{\Rank}{\mathop{\bf Rank}}
\newcommand{\Card}{\mathop{\bf Card}}
\newcommand{\Tr}{\mathop{\bf Tr}}
\newcommand{\diag}{\mathop{\bf diag}}

\newcommand{\Expect}{\mathop{\bf E{}}}

\newcommand{\dsp}{\displaystyle}

\newtheorem{theorem}{Theorem}

\newtheorem{remark}[theorem]{Remark}

\newcounter{exno}

{\begin{quote}\begin{small}\textbf{Proof.}}%
{\end{small}\end{quote}}

{\begin{quote}}{\end{quote}}

\makeatletter
\long\def\@makecaption#1#2{
   \vskip 9pt
   \begin{small}
   \setbox\@tempboxa\hbox{{\bf #1:} #2}
   \ifdim \wd\@tempboxa > 5.5in
        \begin{center}
        \begin{minipage}[t]{5.5in}
        \addtolength{\baselineskip}{-0.95pt}
        {\bf #1:} #2 \par
        \addtolength{\baselineskip}{0.95pt}
        \end{minipage}
        \end{center}
   \else
    \hbox to\hsize{\hfil\box\@tempboxa\hfil}
   \fi
   \end{small}\par
}
\makeatother

\newcounter{oursection}

\newcounter{lecture}

\newtheorem{corollary}{Corollary}
\newtheorem{assumption}{Assumption}

\title{On the Quality of a Semidefinite Programming Bound for Sparse Principal Component Analysis} 

\author{Laurent El Ghaoui \\ EECS Department, UC Berkeley \\ {\tt elghaoui@eecs.berkeley.edu}}
\date{\today}
\begin{document}
\bibliographystyle{plain}

\maketitle
\begin{abstract}
We examine the problem of approximating a positive, semidefinite matrix $\Sigma$ by a dyad $xx^T$, with a penalty on the cardinality of the vector $x$.  This problem arises in sparse principal component analysis, where a decomposition of $\Sigma$ involving sparse factors is sought.  We express this hard, combinatorial problem as a maximum eigenvalue problem, in which we seek to maximize, over a box, the largest eigenvalue of a symmetric matrix that is linear in the variables.  This representation allows to use the techniques of robust optimization, to derive a bound based on semidefinite programming.  The quality of the bound is investigated using a technique inspired by Nemirovski and Ben-Tal (2002).
\end{abstract}

\tableofcontents

\newpage

\subsection*{Notation} The notation $\ones$ denotes the vector of ones (with size inferred from context), while $\Card(x)$ denotes the cardinality of a vector $x$ (number of non-zero elements), and $D(x)$ the diagonal matrix with the elements of $x$ on its diagonal.  We denote by $e_i$ the unit vectors of $\mathbb{R}^n$.  For a $n \times n$ matrix $X$, $X \succeq 0$ means $X$ is symmetric and positive semi-definite.  The notation $B_+$, for a symmetric matrix $B$, denotes the matrix obtained from $B$ by replacing negative eigenvalues by $0$. The notation has precedence over the trace operator, so that $\Tr B_+$ denotes the sum of positive eigenvalues of $B$ if any, and $0$ otherwise. Throughout, the symbol $\Expect$ refers to expectations taken with respect to the normal Gaussian distribution of dimension inferred from context.  Finally, the support of a vector $x$ is defined to be the set of indices corresponding to its non-zero elements.

\section{Introduction}
\label{s:intro}
Given a non-zero $n \times n$ positive semi-definite symmetric matrix $\Sigma$ and a scalar $\rho>0$, we consider the {\em  cardinality-penalized variational problem}
\begin{equation}\label{eq:orig-pb}
    \phi(\rho) := \max_x \: x^T \Sigma x - \rho \Card(x) ~:~ \|x\|_2 = 1 .
\end{equation}
This problem is equivalent to solving the sparse rank-one approximation problem 
\[
\min_z \: \|\Sigma - zz^T \|_F^2 - \rho \Card(z),
\]
which arises in the {\em sparse PCA} problem \cite{ZHT:04,AEJL:04}, where a ``decomposition'' of $\Sigma$ into sparse factors is sought.  We refer to \cite{AEJL:04} for a motivation of the sparse PCA problem, and an overview of its many applications.

In the paper \cite{AEJL:04}, the authors have developed the ``direct sparse PCA'' approach, which leads to the following convex relaxation for the problem (\ref{eq:orig-pb}):
\[
\max_X \: \Tr X \Sigma- \rho \|X\|_1 ~:~ X \succeq 0, \;\; \Tr X = 1.
\]
The above problem is amenable to both general-purpose semidefinite programming (SDP) interior-point codes, and more recent first-order algorithms such as Nesterov's smooth minimization technique \cite{Nes:05}.  Unfortunately, the quality of the relaxation seems to be hard to analyze at present.

In this paper, we introduce two new representations of the problem, and a new SDP bound, based on robust optimization ideas \cite{BeN:02}.  Our main goal is to use the new representations of the problem to analyze the quality of the corresponding bound.  

The paper is organized as follows.  Section \ref{s:equ-vs-ineq} develops some preliminary results allowing to restrict our attention to the case when $\rho < \max_i \Sigma_{ii}$. Section~\ref{s:new-rep} then proposes two new representations for $\phi(\rho)$, one based on largest eigenvalue maximization, and the other on a thresholded version of the Rayleigh quotient.  In section \ref{s:relax-SDP}, we derive an SDP-based upper bound on $\phi(\rho)$, and in section \ref{s:quality}, we analyze its quality: as a function of the penalty parameter $\rho$ first, then in terms of structural conditions on matrix $\Sigma$.

It will be helpful to describe $\Sigma$ in terms of the Cholesky factorization $\Sigma = A^TA$, where $A =[a_1 \ldots a_n]$, with $a_i \in \reals^m$, $i=1,\ldots,n$, where $m = \Rank(\Sigma)$.  Further, we will assume, without loss of generality, that the diagonal of $\Sigma$ is ordered, and none of the diagonal elements is zero, so that $\Sigma_{11} \ge \ldots \ge \Sigma_{nn}>0$.  Finally, we define the set ${\cal I}(\rho) := \{i \::\: \Sigma_{ii} > \rho\}$, and let $n(\rho) := \Card {\cal I}(\rho)$.  

\section{Equality vs.\ Inequality Models}
\label{s:equ-vs-ineq}

In the sequel we will develop SDP bounds for the related quantity
\begin{equation}\label{eq:orig-pb-ineq}
\overline{\phi}(\rho) := \max_x \: x^T \Sigma x - \rho \Card(x) ~:~ \|x\|_2 \le 1 .
\end{equation}
The following theorem says that when $\rho <\Sigma_{11}$, the two quantities $\phi(\rho)$, $\overline{\phi}(\rho)$ are positive and equal; otherwise, both $\phi(\rho)$ and $\overline{\phi}(\rho)$ have trivial solutions.

\begin{theorem}\label{th:ineq-vs-eq}
If $\rho < \Sigma_{11}$, we have $\phi(\rho) = \overline{\phi}(\rho) >0$, and the optimal sets of problems (\ref{eq:orig-pb}) and (\ref{eq:orig-pb-ineq}) are the same. Conversely, if $\rho \ge \Sigma_{11}$, we have $\overline{\phi}(\rho) = 0 \ge \phi(\rho) = \Sigma_{11} - \rho$, and a corresponding optimal vector for $\phi(\rho)$ (resp.\ $\overline{\phi}(\rho)$)  is $x = e_1$, the first basis vector in $\mathbb{R}^n$ (resp.\ $x=0$).
\end{theorem}

\noindent
{\bf Proof:}  If $\rho < \Sigma_{11}$, then the choice $x=e_1$ in (\ref{eq:orig-pb-ineq}) implies $\overline{\phi}(\rho)>0$, which in turn implies that an optimal solution $x^\ast$ for (\ref{eq:orig-pb-ineq}) is not zero.  Since the $\Card$ function is scale-invariant, it is easy to show that without loss of generality, we can assume that $x^\ast$ has $l_2$-norm equal to one, which then results in $\phi(\rho) = \overline{\phi}(\rho)>0$.

Let us now turn to the case when $\rho \ge \Sigma_{11}$.  We develop an expression for $\phi(\rho)$ as follows.
First observe that, since $\Sigma \succeq 0$,
\[
\max_{\|x\|_1=1} \: x^T\Sigma x = \Sigma_{11},
\]
which implies that, for every $x$,
\begin{equation}\label{eq:l1-vs-weighted-l2}
    \Sigma_{11} \|x\|_1^2 \ge x^T\Sigma x.
\end{equation}
Now let $t \ge 0$.  The condition $\phi(\rho) \le -t$ holds if and only if
\[ 
\forall \: x, \; \|x\|_2 = 1 ~:~ \rho \Card(x) \ge t + x^T\Sigma x.
\] 
Specializing the above condition to $x=e_1$, we obtain that $\phi(\rho) \le -t$ implies $\rho \ge \Sigma_{11}+t$.  Conversely, assume that $\rho \ge \Sigma_{11}+t$.  Using (\ref{eq:l1-vs-weighted-l2}), we have for every $x$, $\|x\|_2 = 1$:
\[
\rho \Card(x) \ge \rho \|x\|_1^2 \ge (\Sigma_{11} + t) \|x\|_1^2 \ge x^T\Sigma x + t,
\]
where we have used the fact that $\|x\|_1 \ge 1$ whenever $\|x\|_2 = 1$.  Thus we have obtained that 
$\phi(\rho) \le -t$ with $t \ge 0$ if and only if $\rho \ge \Sigma_{11} +t$, which means that $\phi(\rho) = \Sigma_{11}-\rho$ whenever $\rho \ge \Sigma_{11}$.

Finally, let us prove that $\overline{\phi}(\rho)= 0$ when $\rho \ge \Sigma_{11}$. For every $x \ne 0$ such that $\|x\|_2 \le 1$, we have
\[
\rho \Card(x) \ge \dsp\frac{\rho}{\|x\|_2^2} \|x\|_1^2 \ge \rho \|x\|_1^2 \ge x^T\Sigma x,
\]
which shows that $\overline{\phi}(\rho) \le 0$, and concludes our proof. $\blacksquare$

\bigskip

In the sequel, we will make the following assumption.

\begin{assumption}\label{ass:rho}
We assume that $\rho < \Sigma_{11}$, that is, the set ${\cal I}(\rho):= \{i \::\: \Sigma_{ii} > \rho\}$ is not empty. \end{assumption}

\section{New Representations}
\label{s:new-rep}
\subsection{Largest eigenvalue maximization}  
The following theorem shows that the problem of computing $\phi(\rho)$ can be expressed as a eigenvalue maximization problem, where the sparsity pattern is the decision variable.  

\begin{theorem} \label{th:phi-rep-maxev}
For $\rho \in [0,\Sigma_{11}[$, $\phi(\rho)$ can be expressed as the maximum eigenvalue problem
\begin{equation}\label{eq:phi-rob-B}
    \phi(\rho) = \max_{u \in [0,1]^n} \: \lambda_{\rm max} \left(\sum_{i=1}^n u_i B_i \right) ,
\end{equation}
where $B_i := a_ia_i^T - \rho \cdot I_m$, $i=1,\ldots,n$.  

An optimal solution to the original problem (\ref{eq:orig-pb}) is obtained from a sparsity pattern vector $u$ that is optimal for (\ref{eq:phi-rob-B}), by finding an eigenvector $y$ corresponding to the largest eigenvalue of $D(u)\Sigma D(u)$, and setting $x = D(u)y/\|D(u)y\|_2$, where $D(u) := \diag(u)$.
\end{theorem}

\noindent
{\bf Proof.} Since $\rho < \Sigma_{11}$, the result of Theorem \ref{th:ineq-vs-eq} implies that $\phi(\rho)$ is equal to $\overline{\phi}(\rho)$ defined in in (\ref{eq:orig-pb-ineq}).  Let us now prove that $\phi(\rho) = \tilde{\phi}(\rho)$, where
\begin{equation}\label{eq:phi-uy-rep}
    \tilde{\phi}(\rho) : = \max_{u \in \{0,1\}^n} \: \max_{y^Ty \le 1} \: y^TD(u) \Sigma D(u) y - \rho \cdot \ones^Tu .
\end{equation}

To prove this intermediate result, first note that if $x$ is optimal for $\phi(\rho)$, that is, for (\ref{eq:orig-pb-ineq}), then we can set $u_i = 1$ if $x_i\ne 0$, $u_i=0$ otherwise, so that $\Card(x) = \ones^Tu$; then, we set $y = x$ and obtain that the pair $(u,y)$ is feasible for $\tilde{\phi}(\rho)$, and achieves the objective value $\phi(\rho)$, hence $\phi(\rho) \le \tilde{\phi}(\rho)$.  Conversely, if $(u,y)$ is optimal for $\tilde{\phi}(\rho)$, then $x = D(u)y$ is feasible for $\phi(\rho)$ (as expressed in (\ref{eq:orig-pb-ineq})), and satisfies $\Card(x) \le \Card(u) = \ones^Tu$, thus
\[
\tilde{\phi}(\rho) = y^TD(u)\Sigma D(u) y - \rho \ones^T u \le x^T\Sigma x - \rho \Card(x) \le \phi(\rho),
\]
This concludes the proof that $\phi(\rho) = \tilde{\phi}(\rho)$. 

We proceed by eliminating $y$ from (\ref{eq:phi-uy-rep}), as follows:
\begin{eqnarray*}
  \phi(\rho) &=& \max_{u \in \{0,1\}^n} \: \lambda_{\rm max} (D(u) \Sigma D(u)) - \rho \cdot \ones^Tu  \\  
    &=& \max_{u \in \{0,1\}^n} \: \lambda_{\rm max} (D(u)A^TAD(u)) - \rho \cdot \ones^Tu  \\  
  &=& \max_{u \in \{0,1\}^n} \: \lambda_{\rm max} (AD(u)A^T) - \rho \cdot \ones^Tu  \\
  &=& \max_{u \in \{0,1\}^n} \: \lambda_{\rm max} (\sum_{i=1}^n u_i a_ia_i^T) - \rho \cdot \ones^Tu,
\end{eqnarray*}  
in virtue of $\Sigma = A^TA$, and $D(u)^2 = D(u)$ for every feasible $u$. Invoking the convexity of the largest eigenvalue function, we can replace the set $\{0,1\}^n$ by $[0,1]^n$ in the above expression, and obtain (\ref{eq:phi-rob-B}). $\blacksquare$

\subsection{Thresholded Rayleigh quotient}
\label{ss:alternate-exp} 
The following theorem shows that $\phi(\rho)$ can be expressed as a maximal ``thresholded Rayleigh quotient'', which  for $\rho=0$ reduces to the ordinary Rayleigh quotient.  

\begin{theorem} \label{th:th-rayleigh}
For $\rho \in [0,\Sigma_{11}[$, we have
\begin{eqnarray}
   \phi(\rho) &=&  \max_{\xi^T\xi = 1} \: \sum_{i=1}^n ( (a_i^T\xi)^2 - \rho)_+ , \label{eq:alt-exp-phi} \\
    &=& \max_{X} \: \sum_{i=1}^n ( a_i^TX a_i - \rho)_+  ~:~ X \succeq 0, \;\; \Tr X = 1.  \label{eq:phi-def-dual}
\end{eqnarray}
An optimal solution $x$ for (\ref{eq:orig-pb}) is obtained from an optimal solution $\xi$ to problem (\ref{eq:alt-exp-phi}) by setting $u_i = 1$ if $(a_i^T\xi)^2 > \rho$, $u_i = 0$ otherwise; then, finding an eigenvector $y$ corresponding to the largest eigenvalue of $D(u)\Sigma D(u)$, and setting $x = D(u)y/\|D(u)y\|_2$.
\end{theorem}

\noindent 
{\bf Proof:}  From the expression (\ref{eq:phi-rob-B}), we derive
\begin{eqnarray}
\phi(\rho) &=& \max_{u \in [0,1]^n} \: \max_{\xi^T\xi \le 1} \: \xi^T \left(\sum_{i=1}^n u_i a_ia_i^T\right) \xi - \rho \cdot \ones^Tu \nonumber \\
&=& \max_{\xi^T\xi \le 1} \: \sum_{i=1}^n ( (a_i^T\xi)^2 - \rho \xi^T\xi)_+ \nonumber \\
&=& \max_{\xi^T\xi = 1} \: \sum_{i=1}^n ( (a_i^T\xi)^2 - \rho)_+ , \label{eq:alt-exp-phi-bis}
\end{eqnarray}  
where the last equality derives from the fact that $\phi(\rho)>0$ (which is in turn the consequence of our assumption that $\Sigma_{11} = \max_i a_i^Ta_i > \rho$). Finally, the equivalence between (\ref{eq:alt-exp-phi}) and (\ref{eq:phi-def-dual}) stems from convexity of the objective function in problem (\ref{eq:phi-def-dual}), which implies that without loss of generality, we can impose $X$ to be of rank one in (\ref{eq:phi-def-dual}). $\blacksquare$

\bigskip

The following corollary shows that we can safely remove columns and rows in $\Sigma$ that have variance below the threshold $\rho$.

\begin{corollary}\label{cor:support-soln}
Without loss of generality, we can assume that every optimal solution to the original problem (\ref{eq:orig-pb}) has a support included in the set ${\cal I}(\rho) := \{i \::\: \Sigma_{ii} > \rho\}$.   Thus, if $\Sigma_{ii} \le \rho$, the corresponding column and row can be safely removed from $\Sigma$.
\end{corollary}
\noindent 
{\bf Proof:} 
This is a direct implication of the fact that for every $i$, if $\rho \ge a_i^Ta_i$, then we have $(a_i^T\xi)^2 \le \rho $ for every $\xi$ such that $\xi^T\xi = 1$.  Hence, the corresponding term does not appear in the sum in (\ref{eq:alt-exp-phi-bis}). $\blacksquare$

\subsection{Exact solutions in some special cases}
Theorems \ref{th:phi-rep-maxev} and \ref{th:th-rayleigh} allows to solve exactly the problem in some special cases.

First, Theorem \ref{th:phi-rep-maxev} can be invoked when $\Sigma$ is diagonal, in which case the optimal vector $x$ turns out to be simply the first unit vector, $e_1$.  

Next, consider the case when the matrix $\Sigma$ has rank one, that is, $m=1$.  Then, the $a_i$'s are scalars, and the representation given in Theorem \ref{th:th-rayleigh} yields
\[
\phi(\rho) = \max_{\xi^2 = 1} \: \sum_{i=1}^n ( (a_i \xi)^2 - \rho)_+  = \sum_{i=1}^n ( a_i^2 - \rho)_+ .
\]
A corresponding optimal solution for $\phi(\rho)$ is obtained by setting $u_i = 1$ if $\rho < a_i^2$, $u_i=0$ otherwise, and then setting $x = \tilde{a}/\|\tilde{a}\|_2$, with $\tilde{a}$ obtained from $a$ by thresholding $a$ with absolute level $\sqrt{\rho}$.  In the sequel, we assume that $m>1$.

A similar result holds when $\Sigma$ has the form $\Sigma = I + aa^T$, when $a$ is a given $n$-vector, since then the problem trivially reduces to the rank-one case.

\section{SDP relaxation}  
\label{s:relax-SDP}
A relaxation inspired by \cite{BeN:02} is given by the following theorem.
\begin{theorem} \label{th:sdp-bnd}
For every $\rho \in [0,\Sigma_{11}[$, we have $\phi(\rho) \le \psi(\rho)$, where $\psi(\rho)$ is the solution to a semidefinite program:
\begin{equation}\label{eq:sdp-relax}
\psi(\rho) := \min_{(Y_i)_{i=1}^n} \: \lambda_{\rm max} \left( \sum_{i=1}^n Y_i \right) ~:~ Y_i \succeq B_i , \;\; Y_i \succeq 0, \;\;i=1,\ldots,n .
\end{equation}
The problem can be represented in dual form, as the convex problem
\begin{equation}\label{eq:sdp-dual-X}
\psi(\rho)= \max_X \: \sum_{i=1}^n \Tr \left( X^{1/2} a_ia_i^T X^{1/2} - \rho X \right)_+ ~:~ X \succeq 0, \;\; \Tr X = 1 .
\end{equation}  
\end{theorem}

\noindent 
{\bf Proof:} If $(Y_i)_{i=1}^n$ is feasible for the above SDP, then for every $\xi \in \reals^m$, $\xi^T\xi\le 1$, and $u \in [0,1]^n$, we have
\[
\xi^T \left( \sum_{i=1}^n u_iB_i \right) \xi \leq \sum_{i=1}^n (\xi^T B_i\xi)_+ \le \xi^T \left( \sum_{i=1}^n Y_i \right) \xi \le \psi(\rho), 
\]
which proves $\phi(\rho) \le \psi(\rho)$.  The dual of the SDP (\ref{eq:sdp-relax}) is given by 
\begin{equation}
\psi(\rho) = \max_{X,(P_i)_{i=1}^n} \: \sum_{i=1}^n \langle P_i,B_i \rangle ~:~ X \succeq P_i \succeq 0 , \;\; i=1,\ldots,n , \;\; \Tr X = 1 . 
\end{equation} 
Using the fact that, for any symmetric matrix $B$, and positive semi-definite matrix $X$,
\[
\max_P \: \left\{ \langle P,B \rangle ~:~ X \succeq P \succeq 0 \right\} = \Tr \left( X^{1/2} B X^{1/2} \right)_+ ,
\]
allows to represent the dual problem in the form (\ref{eq:sdp-dual-X}). Note that the convexity of the representation (\ref{eq:sdp-dual-X}) is not immediately obvious. $\blacksquare$

\bigskip

A few comments are in order.

The fact that $\phi(\rho) \le \psi(\rho)$ can also be inferred directly from the dual expression (\ref{eq:sdp-dual-X}): we have, by convexity, and using the representation (\ref{eq:phi-def-dual}) for $\phi(\rho)$, 
\[
\psi(\rho) \ge \max_X \: \left\{ \sum_{i=1}^n \left( a_i^TXa_i - \rho\right)_+ ~:~ X \succeq 0, \;\; \Tr X = 1 \right\} = \phi(\rho) .
\]
From the representation (\ref{eq:sdp-dual-X}) and this, we obtain that if the rank $k$ of $X$ at the optimum of the dual problem (\ref{eq:sdp-dual-X}) is one, then our relaxation is exact: $\phi(\rho) = \psi(\rho)$.  

In fact, problem (\ref{eq:sdp-dual-X}) can be obtained as a rank relaxation of the following exact representation of $\phi$:
\[
\phi = \max_X \: \sum_{i=1}^n \Tr \left( X^{1/2} a_ia_i^T X^{1/2} - \rho X \right)_+ ~:~ X \succeq 0, \;\; \Tr X = 1 , \;\; \Rank(X) = 1.
\]
In contrast, applying a direct rank relaxation to problem (\ref{eq:alt-exp-phi}) (that is, writing the problem in terms of letting $X = \xi\xi^T$ and dropping the rank constraint on $X$) would be useless: it would yield (\ref{eq:phi-def-dual}), which is $\phi(\rho)$ itself.

Finally, note that our relaxation shares the property of the exact formulation (\ref{eq:alt-exp-phi}) observed in Corollary \ref{cor:support-soln}, that indices $i$ such that $\rho \ge \Sigma_{ii}$ can be simply ignored, since then $B_i \preceq 0$.

\section{Quality of the SDP relaxation} 
\label{s:quality} 
In this section, we seek to estimate a lower bound on the {\em quality} of the SDP relaxation, which we define to be a scalar $\theta \in [0,1]$ such that 
\begin{equation}\label{eq:def-rel-approx-error}
\theta \psi(\rho) \le \phi(\rho) \le \psi(\rho).
\end{equation}
Thus, $(1-\theta)/\theta$ is a upper bound on the relative approximation error, $(\psi(\rho)-\phi(\rho))/\phi(\rho)$.

\subsection{Quality estimate as a function of the penalty parameter}
Our first result gives a bound on the relaxation quality conditional on a bound on $\rho$.
We begin by making the following assumption:
\begin{assumption}\label{ass:rho-min}
We assume that $0 < \rho < \dsp\min_{1\le i \le n} \Sigma_{ii} = \Sigma_{nn}$, and $m = \Rank(\Sigma) >1$
.
\end{assumption}
From the result of Corollary~\ref{cor:support-soln}, we can always reduce the problem so that the above assumption holds, by removing appropriate columns and rows of $\Sigma$ if necessary.

\begin{theorem} \label{th:quality-as-fcn-penalty}
With assumption~\ref{ass:rho-min} in force, for every value of the penalty parameter $\rho \in [0,\Sigma_{nn}[$, and for every $\gamma \ge 0$ such that
\begin{equation}\label{eq:rho-u-cond}
    \rho \le \frac{\gamma}{n+\gamma} \Sigma_{11} ,
\end{equation}
the bound (\ref{eq:def-rel-approx-error}) holds with $\theta$ set to $\theta_{m}(\gamma)$, where for $m>1$ and $\gamma \ge 0$, we define 
\begin{equation}\label{eq:def-theta-u}
\begin{array}{rcl}
    \theta_m(\gamma) &:=& 
    \Expect \left( \xi_1^2 - \dsp\frac{\gamma}{m-1} \cdot \dsp\sum_{j=2}^{m} \xi_{j}^2 \right)_+ ,
\end{array}
\end{equation}
which can be computed by the formula
\begin{equation}\label{eq:def-theta-rep}
\begin{array}{rcl}    
\theta_m(\gamma) &=& \frac{\dsp\int_0^{\pi/2} \left(\cos^2(t) -\frac{\gamma}{m-1} \sin^2(t)\right)_+ \sin^{m-2}(t)dt}{
\dsp\int_0^{\pi/2} \cos^2(t) \sin^{m-2}(t)dt} .
\end{array}    
\end{equation}
For every $\gamma \ge 0$, the value $\theta_m(\gamma)$ decreases with $m$, and admits the bound
\begin{equation}\label{eq:theta-bnd}
    \theta_m(\gamma) \ge \frac{1}{2} \left( 1 - \gamma + \frac{2}{\pi} \sqrt{1 + \frac{\gamma^2}{m-1} } \right)_+  .
\end{equation}

In particular, if $\rho$ satisfies (\ref{eq:rho-u-cond}) with $\gamma = 1$, that is, $\rho \le \Sigma_{11}/(n+1)$, then bound (\ref{eq:def-rel-approx-error}) holds with $\theta \ge 1/\pi$.
\end{theorem}

Before we prove the theorem, let us make a few comments.  

First, as will be apparent from the proof, the value of $m$ can be safely replaced by the rank $k$ of an optimal solution to the SDP (\ref{eq:sdp-dual-X}).  This can only improve the quality estimate, as $k \le m$ and $\theta_m(\gamma)$ is a decreasing function of $m$ for every $\gamma \ge 0$.

\begin{figure}[thp]
\begin{center}
\psfrag{t}{}
\psfrag{x}[t][][1][0]{$\gamma$}
\psfrag{y}[r][][1][270]{$\theta_m(\gamma)$}
\psfrag{k=2}{\tiny $m=2$}
\psfrag{k=3}{\tiny $m=3$}
\psfrag{k=10}{\tiny $m=10$}
\psfrag{k=100}{\tiny $m=100$}
\begin{tabular}{c}
\includegraphics[width=0.55 \linewidth,height=.3 \textheight]{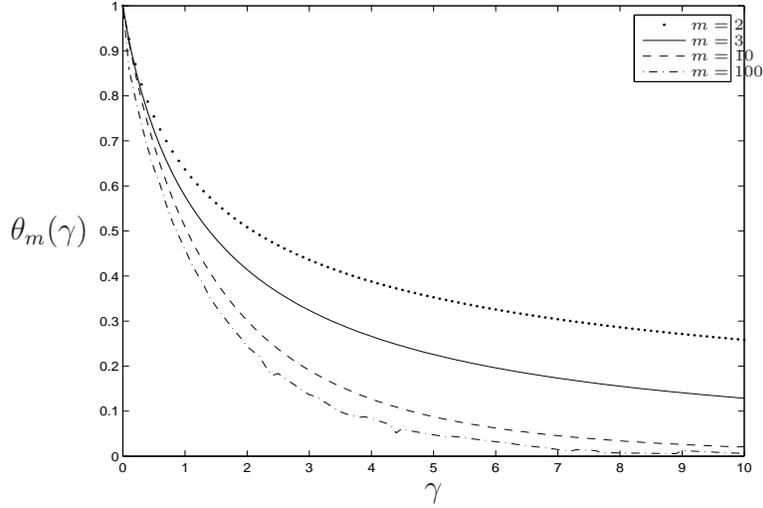} 
\end{tabular}
\caption{Plot of function $\theta_m(\gamma)$, as defined in (\ref{eq:def-theta-rep}), for various values of $m$.}
\label{fig:theta-vs-gamma}
\end{center}
\end{figure}

Second, the smaller $m$ is, and the larger $\gamma$ is, the smaller the corresponding quality estimate.  However, a small value for $\gamma$ does not allow for a large range of $\rho$ values via (\ref{eq:rho-u-cond}), and this effect is becomes more pronounced as $n$ grows. The theorem presents the result in such a way that the respective contributions of $m,n$ to the deterioration of the quality estimate are separated.  A plot of the function $\theta_m$ for various values of $m$ is shown in Figure~\ref{fig:theta-vs-gamma}.  

Third, the theorem allows to plot the predicted quality estimate $\theta$ as a function of the penalty parameter, in the interval $[0,\Sigma_{nn}[$.  Leveraging these results to the entire range $[0,\Sigma_{11}[$ will be straightforward, but will require us to be careful about the sizes $n$ and $m$, as they change as $\rho$ crosses the values $\Sigma_{n-1,n-1},\ldots,\Sigma_{11}$, in view of Corollary \ref{cor:support-soln}.  We formalize the argument in Corollary \ref{cor:quality-fcn-of-rho}.

Finally, the theorem allows to derive conditions on the structure of $\Sigma$ that guarantee a prescribed value of the quality.  We describe such a condition in Corollary \ref{cor:quality-structural}.

\bigskip

\noindent{\bf Proof of theorem \ref{th:quality-as-fcn-penalty}:} The approach we use in our proof is inspired by that of Theorem 2.1 in \cite{BeN:02}.  Let $X \succeq 0$, $\Tr X = 1$, be optimal for the upper bound $\psi(\rho)$ in dual form (\ref{eq:sdp-dual-X}), so that 
\[
\psi(\rho) = \sum_{i=1}^n \Tr (B_i(X)_+) ,
\]
where $B_i(X) := X^{1/2}B_i X^{1/2}$.  Let $k = \Rank(X)$.  We have seen that if $k=1$, then our relaxation is exact: $\phi(\rho) = \psi(\rho)$.  If the rest of the proof, we will assume that $k>1$.  We thus have $1 < k \le m = \Rank (\Sigma) \le n$.

Assume that we find a scalar $\theta \in [0,1]$ such that:
\begin{equation}\label{eq:bound-in-expect}
    \Expect \sum_{i=1}^n (\xi^TB_i(X)\xi)_+ > (\theta \psi(\rho)) \cdot \Expect (\xi^TX\xi) ,
\end{equation}
where $\xi$ follows the normal distribution in $\mathbb{R}^m$. The bound above implies that there exist a non-zero $\xi \in \reals^m$ such that
\[
\sum_{i=1}^n  (\xi^TB_i(X)\xi)_+ > (\theta \psi(\rho)) \cdot  (\xi^TX \xi) .
\]
Thus, with $u_i = 1$ if $\xi^TB_i(X)\xi >0$, $u_i=0$ otherwise, we obtain that there exist a non-zero $\xi \in \reals^m$ and $u  \in [0,1]^n$ such that
\[
\xi^T \left( \sum_{i=1}^n  u_i B_i(X) \right) \xi > (\theta \psi(\rho)) \cdot  (\xi^TX \xi) .
\]
With $z = X^{1/2}\xi$:
\[
z^T \left( \sum_{i=1}^n  u_i B_i \right) z > (\theta \psi(\rho)) \cdot (z^Tz) .
\]
The above implies that $z \ne 0$, so we conclude that there exist  $u \in [0,1]^n$ such that
\[
\lambda_{\rm max}\left( \sum_{i=1}^n  u_i B_i \right) > \theta \psi(\rho),
\]
from which we obtain the quality estimate $\theta \psi(\rho) \le \phi(\rho) \le \psi(\rho)$.  By a continuity argument, this result still holds if (\ref{eq:bound-in-expect}) is satisfied, but not strictly.  The rest of the proof is dedicated to finding a scalar $\theta$ such that the bound (\ref{eq:bound-in-expect}) holds.

Fix $i \in \{1,\ldots,n\}$.  It is easy to show that $B_i(X)$ has exactly one positive eigenvalue $\alpha_i$, since assumption \ref{ass:rho-min} holds. Thus $\alpha_i = \Tr B_i(X)_+$.  Since $B_i(X) \preceq X^{1/2}a_ia_i^TX^{1/2}$, we have $\alpha_i = \lambda_{\rm max}(B_i(X)) \le a_i^TXa_i$. Further, $B_i(X)$ has exactly rank $k=\Rank(X)$.  Denote by $(-\beta_j^i)_{j=1}^{k-1}$ the negative eigenvalues of $B_i(X)$. We then have
\[
\sum_{j=1}^{k-1} \beta_j^i = \Tr B_i(X)_+ - \Tr B_i(X) = \alpha_i - (a_i^TXa_i-\rho) \le \rho.
\]

Now let $\xi$ follow the normal distribution in $\mathbb{R}^m$, ${\cal N}(0,I_m)$.  By rotational invariance of the normal distribution, we have:
\begin{eqnarray*}
\Expect (\xi^T B_i(X) \xi)_+ &=& \Expect \left( \alpha_i \xi_1^2 - \sum_{j=1}^{k-1} \beta_j^i \xi_{j+1}^2 \right)_+ .
\end{eqnarray*}
Thus,
\begin{eqnarray}
\Expect (\xi^T B_i(X) \xi)_+ &\ge& \dsp\min_{\beta \in \mathbb{R}^k} \: \left\{ \Expect \left( \alpha_i \xi_1^2 - \sum_{j=1}^{k-1} \beta_{j} \xi_{j+1}^2 \right)_+ ~:~ \beta \ge 0, \;\; \sum_{j=1}^{k-1} \beta_j \le \rho  \right\} \label{eq:expect-rep-min} \\
&\ge& \dsp\min_{\beta \in \mathbb{R}^m} \: \left\{ \Expect \left( \alpha_i \xi_1^2 - \sum_{j=1}^{m-1} \beta_{j} \xi_{j+1}^2 \right)_+ ~:~ \beta \ge 0, \;\; \sum_{j=1}^{m-1} \beta_j \le \rho  \right\} \label{eq:expect-rep-min-m} \\
&=& \Expect \left( \alpha_i \xi_1^2 - \dsp\frac{\rho}{m-1} \sum_{j=1}^{m-1} \xi_{j+1}^2 \right)_+ , \label{eq:expect-rep-min-sol}
\end{eqnarray}
where we have exploited the convexity and symmetry in problem (\ref{eq:expect-rep-min-m}).  (As claimed in the first remark made after Theorem \ref{th:quality-as-fcn-penalty}, we could safely keep $k$ instead of $m$ in the remaining of the proof.)

Summing over $i$, and in view of $\psi(\rho) = \sum_{i=1}^n \alpha_i$, we get:
\begin{eqnarray}
\Expect \sum_{i=1}^n (\xi^TB_i(X)\xi)_+ 
&\ge& \sum_{i=1}^n \Expect \left( \alpha_i \xi_1^2 - \dsp\frac{\rho}{m-1} \sum_{j=1}^{m-1} \xi_{j+1}^2 \right)_+ \nonumber \mbox{ (by the bound (\ref{eq:expect-rep-min-sol}))} \\
&\ge& \Expect \left( \psi(\rho) \xi_1^2 - \dsp\frac{n \rho}{m-1} \sum_{j=1}^{m-1} \xi_{j+1}^2 \right)_+ \nonumber \mbox{ (by homogeneity and convexity)} \\
&\ge& \theta_m(\gamma) \cdot \psi(\rho), 
\end{eqnarray}
provided $\gamma \ge n \rho / (m-1)\psi(\rho)$. Using the fact that $\psi(\rho) \ge \phi(\rho) \ge \Sigma_{11} - \rho$, we obtain that the bound (\ref{eq:def-rel-approx-error}) holds with $\theta = \theta_m(\gamma)$ whenever (\ref{eq:rho-u-cond}) does, as claimed in the theorem.  The expression (\ref{eq:def-theta-rep}) of the function $\theta_m$ is proved in Appendix \ref{app:theta-rep}, while the bound (\ref{eq:theta-bnd}) is proved in Appendix \ref{app:theta-bound}. Finally, the fact that the function $\theta_m(\gamma)$ decreases with $m>1$ for every $\gamma>0$ is a consequence of the following representation:
\[
\theta_m(\gamma) = \dsp\min_{\beta \in \mathbb{R}^m} \: \left\{ \Expect \left(  \xi_1^2 - \sum_{j=1}^{m-1} \beta_{j} \xi_{j+1}^2 \right)_+ ~:~ \beta \ge 0, \;\; \sum_{j=1}^{m-1} \beta_j \le \gamma  \right\} .
\]
Indeed, adding constraints $\beta_j = 0$ for $j>k$ in the above problem shows that $\theta_m(\gamma) \le \theta_k(\gamma)$ for every $\gamma >0$ and $k \le m$. 
$\blacksquare$

\bigskip

\begin{figure}[thp]
\begin{center}
\psfrag{t}{}
\psfrag{x}[t][][1][0]{$\rho/\Sigma_{11}$}
\psfrag{y}[r][][1][270]{$\vartheta(\rho)$}
\begin{tabular}{ccc}
\includegraphics[width=0.45 \linewidth,height=.3 \textheight]{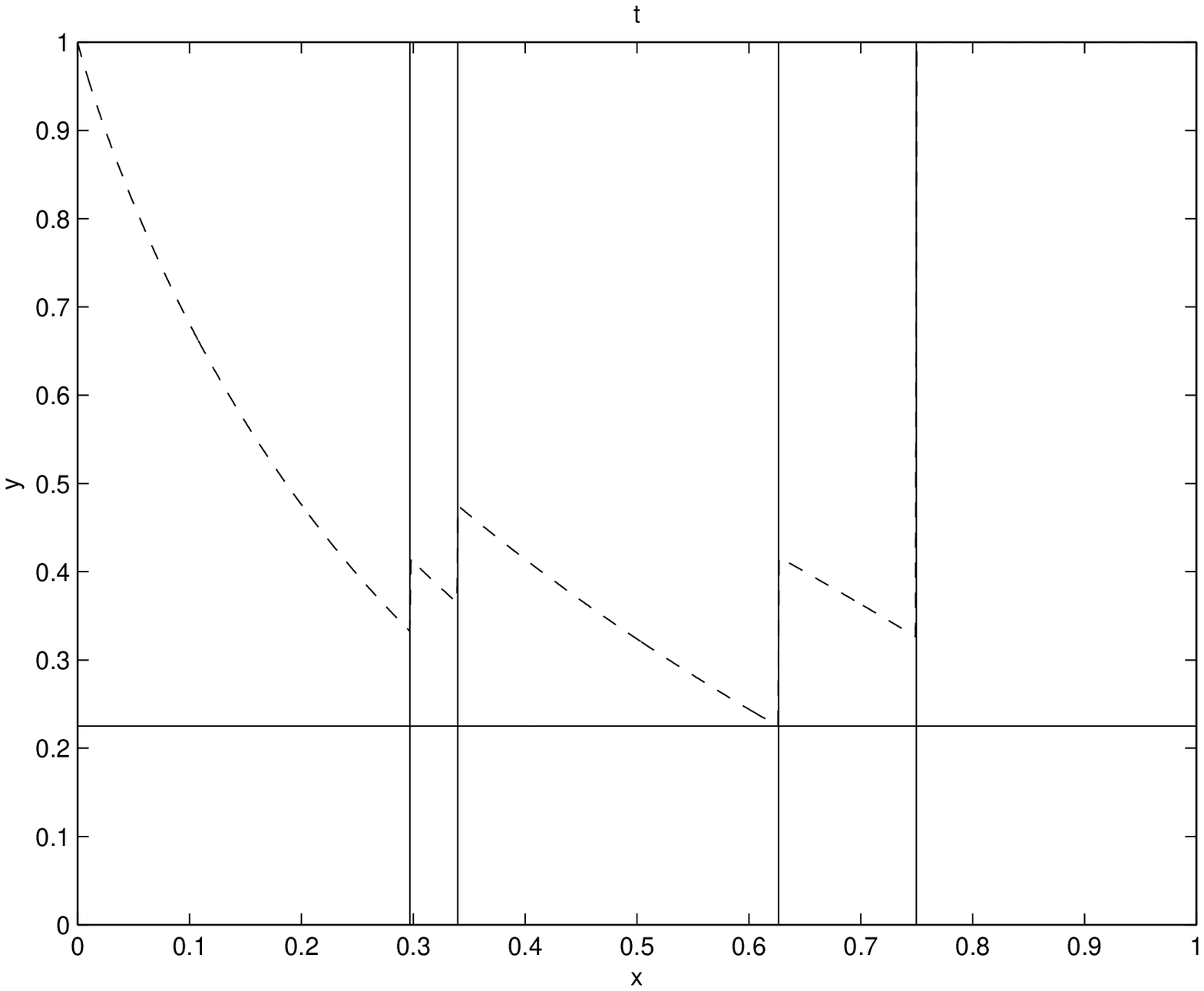} 
& \hspace{.05in} &
\includegraphics[width=0.45 \linewidth,height=.3 \textheight]{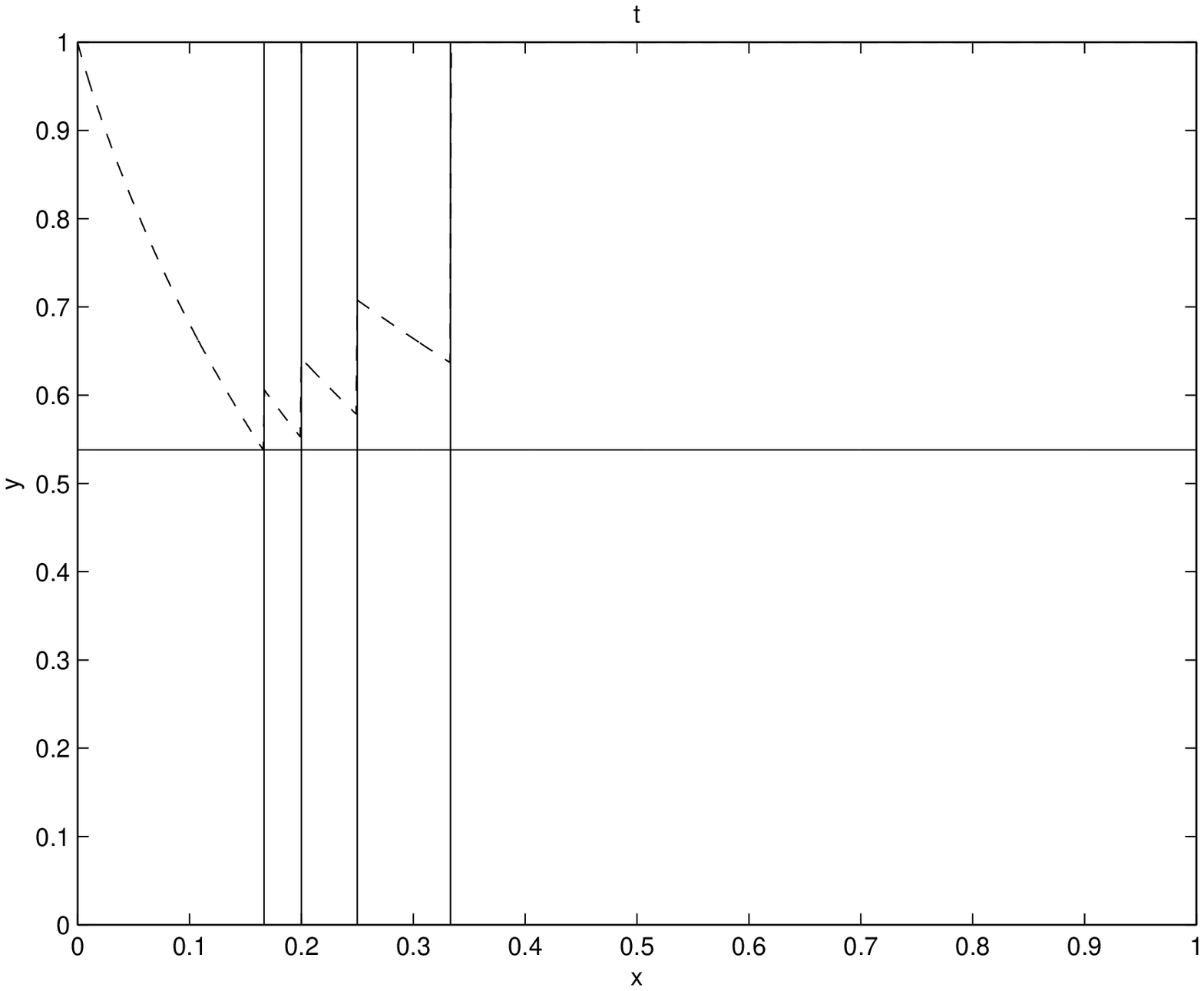} 
\end{tabular}
\caption{Plot of the function $\vartheta(\rho)$ defined in Corollary \ref{cor:quality-fcn-of-rho}, for a specific $5 \times 5$ covariance matrix $\Sigma$.  The left pane corresponds to a random matrix, and the right pane, to a random matrix that satisfies the conditions of Corollary \ref{cor:quality-structural}.}
\label{fig:theta-vs-rho}
\end{center}
\end{figure}

The following corollary allows to plot the quality estimate, as derived from Theorem~\ref{th:quality-as-fcn-penalty}, as a function of $\rho$ across the entire range $[0,\Sigma_{11}[$.  We do not make the assumption~\ref{ass:rho-min} anymore, but do keep assumption \ref{ass:rho}.

\begin{corollary}\label{cor:quality-fcn-of-rho}
Let $\rho \in [0,\Sigma_{11}[$, and define $n(\rho) = \Card\{i \::\: \Sigma_{ii} > \rho\}>0$ and $m(\rho) = \Rank(\Sigma(\rho))$, where $\Sigma(\rho)$ is the $n(\rho) \times n(\rho)$ matrix obtained by removing the last $n-n(\rho)$ rows and columns in $\Sigma$. The bound (\ref{eq:def-rel-approx-error}) holds for $\theta = {\bf \vartheta}(\rho)$, where
\[
\vartheta(\rho) = \left\{ \begin{array}{ll} 
\theta_{m(\rho)}(\gamma(\rho)), \;\; \gamma(\rho) = \dsp\frac{n(\rho)}{m(\rho)-1} \cdot \dsp\frac{\rho}{\Sigma_{11}-\rho} & \mbox{if } m(\rho) >1, \\
1 & \mbox{otherwise.}
\end{array} \right.
\]
\end{corollary}
An example of the resulting plot is shown in Figure~\ref{fig:theta-vs-rho}.

\subsection{Quality estimate based on the structure of $\Sigma$}
\label{ss:quality-structure}

The next result illustrates how to obtain a quality estimate based on structural assumptions on $\Sigma$, requiring that its ordered diagonal decreases fast enough. 

\begin{corollary}\label{cor:quality-structural}
Assume $\Sigma_{11} > \ldots > \Sigma_{nn}$. If $\Sigma_{22} \le \rho < \Sigma_{11}$, 
then the bounds (\ref{eq:def-rel-approx-error}) hold with $\theta =1$, that is, $\phi(\rho) = \psi(\rho)$.  If in addition, we have, for every $h \in \{2,\ldots,n\}$
\begin{equation}\label{eq:quality-cond-struc-x}
   \Sigma_{hh} \le \frac{1}{h+1}  \Sigma_{11}, 
\end{equation}
then, whenever $0 < \rho < \Sigma_{22}$, the bounds (\ref{eq:def-rel-approx-error}) hold with $\theta \ge 1/\pi$.
\end{corollary}

\bigskip

\noindent
{\bf Proof:}  In the case $\rho \in [\Sigma_{22},\Sigma_{11}[$, $n(\rho) = 1$, so that $m(\rho) = 1$, and the bound (\ref{eq:def-rel-approx-error}) holds with $\theta =1$.  Now let $\rho$ be such that $0 < \rho < \Sigma_{22}$. Then there exist $h \in \{2,\ldots,n\}$ such that $\Sigma_{h+1,h+1} \le \rho < \Sigma_{hh}$, with the convention $\Sigma_{n+1,n+1} = 0$.  In  this case, $n(\rho) = \Card\{i \::\: \Sigma_{ii} > \rho\} = h$, so that the sufficient condition (\ref{eq:rho-u-cond}) with $\gamma = 1$ writes
\[
\rho \le \dsp\frac{1}{(h+1)} \Sigma_{11},
\]
which, in view of $\Sigma_{h+1,h+1} \le \rho < \Sigma_{hh}$, holds when (\ref{eq:quality-cond-struc-x}) holds, independent of $\rho$.  Applying the bound (\ref{eq:theta-bnd}) ends the proof. $\blacksquare$

An example corresponding to the situation of Corollary \ref{cor:quality-structural} is shown in Figure~\ref{fig:theta-vs-rho} (left pane).


\bibliography{maxev}

\begin{appendix}
\section{A Formula for $\theta_m$}
\label{app:theta-rep}
Let $\gamma>0$, $m>1$.  Let us prove that the function $\theta_m$ defined in (\ref{eq:def-theta-u}) can be represented as in (\ref{eq:def-theta-rep}). Using the hyperspherical change of variables 
\[
\begin{array}{rcl}
\xi_1 &=& r \cos(\phi_1) \\
\xi_2 &=& r \sin (\phi_1) \cos(\phi_2) \\
&\vdots& \\
\xi_{m-1} &=& r \sin (\phi_1) \ldots \sin(\phi_{m-2}) \cos(\phi_{m-1}) \\
\xi_{m} &=& r \sin (\phi_1) \ldots \sin(\phi_{m-2}) \sin(\phi_{m-1}) ,
\end{array}
\]
with $\phi_j \in [0,\pi]$, $j=1,\ldots,m-2$, $\phi_{m-1} \in [0,2\pi]$, and 
with the corresponding change of measure 
\[
d\xi = r^{m-1} \sin^{m-2}(\phi_1) \ldots \sin(\phi_{m-2}) d\phi_1 \ldots d\phi_{m-1} , 
\]
we obtain
\begin{eqnarray*}
\theta_m(\gamma) &=& (2\pi)^{-{m}/2} \dsp\int_{\mathbb{R}^{m}} \left( \xi_1^2 - \frac{\gamma}{m-1} \sum_{j=2}^{{m}} \xi_{j}^2 \right)_+ e^{-\|\xi\|_2^2/2} d\xi \\
&=&  I_{m} \cdot J_m(\gamma) , 
\end{eqnarray*}
where $I_m$ is some constant, independent of $\gamma$, and
\[
J_m(\gamma) := \int_{0}^\pi  \left(\cos^2(\phi_1) -  \frac{\gamma}{m-1}\sin^2(\phi_1) \right)_+ \sin^{m-2}(\phi_1) d\phi_1. 
\]
Since $\theta_m(0) = 1$, we have $I_{m} = 1/J_{m}(0)$. 
Exploiting symmetry to reduce the integration interval from $[0,\pi]$ to $[0,\pi/2]$, proves the formula (\ref{eq:def-theta-rep}).

\section{A bound on $\theta_m$}
\label{app:theta-bound}
The bound stems from the identity $a_+ = (a+|a|)/2$, valid for every $a \in \mathbb{R}$, and the following result, found in the proof of Theorem 2.1 of \cite{BeN:02}:
\[
\forall \: y \in \mathbb{R}^m ~:~ \Expect \left|\sum_{i=1}^m y_i \xi_i^2 \right| \ge \dsp\frac{2}{\pi} \|y\|_2  .
\]
\end{appendix}

\end{document}